\documentclass[10pt]{article} 
\usepackage{lineno}
\usepackage{amssymb} 
\usepackage{etoolbox}
\usepackage{amsmath} 
\usepackage{amsthm} 
\usepackage{hyperref} 
\usepackage{amsfonts} 
\usepackage{color} 
\usepackage{epsfig} 
\usepackage{graphics,subfigure} 
\usepackage{graphicx} 
\usepackage{float} 
\usepackage{epsf} 
\usepackage{mathrsfs}
\usepackage[numbers,sort&compress]{natbib}
\usepackage{appendix}
\definecolor{myGreen}{rgb}{0.9, 0.99, 0.9}
\newcommand*\linenomathpatch[1]{%
  \cspreto{#1}{\linenomath}%
  \cspreto{#1*}{\linenomath}%
  \csappto{end#1}{\endlinenomath}%
  \csappto{end#1*}{\endlinenomath}%
}

\linenomathpatch{equation}
\linenomathpatch{gather*}
\linenomathpatch{multline}
\linenomathpatch{align}
\linenomathpatch{alignat}
\linenomathpatch{flalign}
\linenomathpatch{cases}

\newtheorem{theorem}{Theorem}[section]
\newtheorem{definition}{Definition}[section]
\newtheorem{lemma}[theorem]{Lemma}
\newtheorem{proposition}[theorem]{Proposition}
\newtheorem{remark}{Remark}[section] 

\allowdisplaybreaks 

\newcommand{\nn}{\nonumber} 

\begin{document}
\begin{titlepage}
\title{\bf Global regularity of 2D Rayleigh-B\'{e}nard equations with logarithmic supercritical dissipation}
\author{ Baoquan Yuan\thanks{Corresponding Author: B. Yuan}\ \ ,\ Xinyuan Xu,\ Changhao Li
       \\ School of Mathematics and Information Science,
       \\ Henan Polytechnic University,  Henan,  454000,  China.\\
      (bqyuan@hpu.edu.cn; xxy11251999@163.com; lch1162528727@163.com)
          }

\date{}
\end{titlepage}
\maketitle
\begin{abstract}

In this paper, we study the global regularity problem for the 2D Rayleigh-B\'{e}nard equations with logarithmic supercritical dissipation. By exploiting a combined quantity of the system, the technique of Littlewood-Paley decomposition and Besov spaces,  and some commutator estimates, we establish the global regularity of a strong solution to this equations in the Sobolev space $H^{s}(\mathbb{R}^{2})$ for $s \ge2$.
\vskip 0.2in
\end{abstract}

\vspace{.2in} {\bf Key words:}\quad  Rayleigh-B\'{e}nard equations; global regularity; logarithmic supercritical dissipation; Sobolev space.

\vspace{.2in} {\bf MSC(2010):} 35Q35; 76D03; 35B65.


\section{Introduction}
	\setcounter{equation}{0}

The aim of this paper is to establish the global regularity of the 2D Rayleigh-B\'{e}nard equations with velocity logarithmic supercritical dissipation
\begin{eqnarray}\label{1.1}
	\left\{\begin {array}{lll}
	\partial_{t}u+ u \cdot \nabla u +\mathcal{L}u + \nabla p = \theta e_{2}, \\
	\partial_{t} \theta + u \cdot \nabla \theta =u \cdot e_{2},\\
	\nabla \cdot u =0,\\
	u(x,0)=u_{0}(x) , \quad \theta(x,0) = \theta_{0}(x),
	\end{array}\right.
	\end{eqnarray}
where $u(x,t)=(u_{1}(x,t),u_{2}(x,t))$ is the velocity field, $\theta = \theta(x,t)$ is a scalar function denoting the temperature and $p$ is the scalar pressure, $e_{2}$ is a unit vector in the $x_{2}$ direction. The forcing term $\theta e_{2}$ describes the acting of the buoyancy force on fluid motion and $u \cdot e_{2}$ models the Rayleigh-B\'{e}nard convection in a heated inviscid fulid, and the operator $\mathcal{L}$ is a Fourier multiplier with symbol $\frac{|\xi|}{g(\xi)}$ and defined as	
 \begin{equation}\label{1.2}
    \hat{\mathcal{L}u}(\xi) = \frac{|\xi|}{g(\xi)}\hat{u}(\xi),
    \end{equation}  	
where $g(\xi)=g(|\xi|)$ is a radial non-decreasing smooth function satisfying the following three conditions:	
\par
\hangindent=3em
\hangafter=1
 (a)  \begin{align}
  &g(\xi) \geq C_{0}>0 \quad for\ all \quad \xi \geq 0; \nn
    \end{align}
\par
\hangindent=3em
\hangafter=1	
 (b) $g$ is of the Mikhlin-H$\rm\ddot{o}$mander type function. That is: there exists a constant $C>0$ such that
    \begin{align}
  &|\partial_{\xi}^{k}g(\xi)| \leq C|\xi|^{-k}|g(\xi)|, \quad k \in \{1,2\}, \quad \forall \ \xi \neq 0;\nn
    \end{align}
    \par
  \hangindent=3em
\hangafter=1	
 (c) $g$ satisfies
    \begin{align}\label{1.3}
   & \lim\limits_{|\xi| \rightarrow \infty} \frac{g(|\xi|)}{|\xi|^{\sigma}}=0 , \quad \forall \ \sigma >0.
   \end{align}

\begin{remark}\label{Remark1.1}
The typical examples of $g$ satisfying the conditions $(a)-(c)$ are
\begin{equation*}
\begin{split}
&g(r)=(\ln(e+r))^{\mu_1}\quad \mbox{for} \ \mu_{1}> 0;\\
&g(r)=\ln(e+r)(\ln(e^2+\ln(1+r)))^{\mu_2}\quad \mbox{for} \ \mu_{2}> 0.
\end{split}
\end{equation*}
\end{remark}\par
The B\'{e}nard problem comes from the convective motions in a heated and incompressible
fluid, which can be used to model the heat convection phenomenon (see, e.g. \cite{AME2015,COR87,TS07,P1968}). According to the classification of the generalized Boussinesq equation (see \cite{QCJZ2014}), we can divide the velocity dissipation $\Lambda^{\alpha}u$ into three cases: subcritical case $\alpha>1$, critical case $\alpha=1$ and supercritical case $\alpha <1$. The global regularity of a strong solution in the subcritical dissipation case is easy, and in 2017, Ye \cite{Z2017} obtained the global well-posedness for the system (\ref{1.1}) with velocity critical dissipation. For the supercritical dissipation case the global regularity of a strong solution is a very difficulty open problem. In this paper, we aim to explore how far one can go beyond the velocity critical dissipation, and we find that the global regularity of a strong solution can still be proved when the critical velocity dissipation $\Lambda u$ reduces to a logarithmic supercritical dissipation $\mathcal{L}u$ as defined in (\ref{1.2}). When the equations $({\ref{1.1}})$ are without the Rayleigh-B\'{e}nard convection term $u\cdot e_{2}$, it reduces to the corresponding  2D Boussinesq equations, which models geophysical flows such as atmospheric fronts and oceanic circulation, and plays an important role in the study of Rayleigh-B\'{e}nard convection.  The Boussinesq equations have recently attracted considerable attention in the community of mathematical fluids (see, e.g. \cite{DCJ2011,AH2020,Z2014,CL2011} and the references therein). In 2024, KC, Regmi, Tao and Wu \cite{DLJD} obtained the global regularity of the Boussinesq equations with only velocity logarithmic supercritical dissipation. For the Rayleigh-B\'{e}nard equation, there is a  Rayleigh-B\'{e}nard convection term $u\cdot e_{2}$ in the temperature $\theta$ equation, which brings about some difficulties. By overcoming these difficulties we establish the global regularity of the Rayleigh-B\'{e}nard equations with velocity logarithmic supercritical dissipation $\mathcal{L} u$. To this end, we state the main result as in the following theorem.
\begin{theorem}\label{thm 1.1}
Consider the system ({\ref{1.1}}) and assume that $\mathcal{L} u$ is defined by ({\ref{1.2}})
with $g(|\xi|)$ obeying (a)-(c). Let the initial data $(u_{0}, \theta_{0})$ be in the class
    \begin{align}
      &u_{0} \in H^{s}(\mathbb{R}^{2}), \quad \theta_{0} \in H^{s}(\mathbb{R}^{2}) \nn
      \end{align}
for $s\ge 2$. Then the system ({\ref{1.1}}) admits a unique global strong solution such that for any $T > 0$,
    \begin{align}
      &u \in C([0,T);H^{s}(\mathbb{R}^{2})) \cap L^{2}([0,T);B_{2,2}^{s+\frac{1}{2},g^{-\frac{1}{2}}}(\mathbb{R}^{2})) ,\quad \theta \in C([0,T);H^{s}(\mathbb{R}^{2})). \nn
      \end{align}
\end{theorem}
\begin{remark}\label{Remark1.2}
$B_{2,2}^{s+\frac{1}{2},g^{-\frac{1}{2}}}(\mathbb{R}^{2})$ denotes a Besov type space, whose  definition is given in the Definition {\ref{def 2.1}}, which is a little weaker than the Besov space $B_{2,2}^{s+\frac{1}{2}}(\mathbb{R}^{2})$.
\end{remark}

\begin{remark}
The global regularity result of a strong solution in Ye \cite{Z2017} is improved by Theorem \ref{thm 1.1}, in which we reduce the critical velocity dissipation $\Lambda u$ to a logarithmic supercritical dissipation $\mathcal{L} u$.
\end{remark}

We now explain the main ideas and the main difficulties encountered in studying the global regularity of solutions of the system ({\ref{1.1}}). A key step to prove the global regularity of solutions to (\ref{1.1}) is to establish a suitable global $\textit{a priori}$ bounds for the solution. To obtain the $H^{1}$ estimate of $u$, we make use the vorticity $\omega=\nabla \times u$ which satisfies
    \begin{align}
      &\partial_{t}\omega + u \cdot \nabla \omega + \mathcal{L}\omega = \partial_{x_1}\theta. \nn
      \end{align}
Since $\theta$ has no viscous derivative, $\partial_{x_1}\theta$ cannot be directly treated by the energy method. In order to overcome this difficulty, we introduce a new combination quantity by using the method in \cite{TSF2010}. More precisely, we set $\mathcal{R}_{g} \triangleq \mathcal{L}^{-1}\partial_{x_1}$ which is similar to a two-dimensional classical Riesz operator and hide the term $\partial_{x_1} \theta$ in the combined quantity
    \begin{equation}\label{1.4}
     G \triangleq \omega-\mathcal{R}_{g}\theta
           \end{equation}
which obeys the following equation
    \begin{equation}\label{1.5}
      \partial_{t}G + u \cdot \nabla G + \mathcal{L}G = [\mathcal{R}_{g},u \cdot \nabla]\theta + \mathcal{R}_{g}u_{2}.
      \end{equation}
Due to the presence of $u \cdot e_{2}$ in the temperature equation, we can not do the estimating $\|\theta\|_{L^\infty}$ directly, but we can obtain the  $L^{p}-$estimate of $\theta$ for $ 2\leq p < 4$. With this estimate at disposal, by a commutator estimate Lemma \ref{lem 2.7}, we can obtain the estimate of $\|G\|_2$, and further the estimate of $\int_{0}^{t}\|\omega\|_{B_{p,\infty}^{0,g^{-1}}}^{2} \mathrm{d}\tau$, and the estimate $\|\theta\|_{L^{\infty}}$ can be obtained. Next, by taking advantage of Littlewood-Paley theory and Besov spaces, we can get the estimate of $\|G\|_{L_{t}^{1}B_{r,1}^{\frac{2}{r}}}$ with $r \in (2,4)$. Furthermore, by utilizing the structure of combined quantity $G$, we can show the estimate of $\|\omega\|_{L_{t}^{1}B_{\infty,1}^{0}}$, which implies $\|\nabla u\|_{L_{t}^{1}L^{\infty}} < \infty$. Finally, we are able to show the global $H^s$ estimates of $u$ and $\theta$ by using commutator estimates and the \textit{a priori} estimates we have obtained.

The rest of this paper is divided into three sections. In Section 2, we recall the Littlewood-Paley theory and the definition of Besov type spaces, and some useful commutator estimates. In Section 3,
we do some $\textit{a priori}$ estimates for strong solutions to the system ({\ref{1.1}}). This section is divided into four subsections. First we do the $L^2-$estimates of $u$ and $\theta$, the $L^2-$estimate of $G$
and the $L^\infty-$ estimate of $\theta$. The estimates of $\|G\|_{L^{r}}$ and $\|G\|_{L^{1}_{t}{B_{r,1}^{\frac{2}{r}}}}$ for $2<r<4$ are carried out in Section 3.2. We further do the estimates of $\|\omega\|_{L^{1}_{t}{B_{\infty,1}^{0}}}$, $\|\nabla u\|_{L_{t}^{1}L^\infty}$ and $\|\nabla \theta\|_{L^\infty}$ in Section 3.3. Finally, we can obtain the global $H^{s}$ estimates of $u$ and $\theta$ in Section 3.4 for $s\ge 2$. In Section 4, by proving the existence and uniqueness of a strong solution to (\ref{1.1}), we complete the proof of Theorem ${\ref{thm 1.1}}$.

	
\section{Preliminaries}
\setcounter{equation}{0}
\vskip .1in

\ \ \ \ In this section we introduce Besov type spaces and several tool lemmas which serve as a preparation for the proofs of our main results.\par
 Now we recall the Littlewood-Paley decomposition and their elementary properties. Related materials can be found in several books and many papers (see e.g \cite{HJR2011,J1998,TSF2011,RM2011}). Let $\mathcal{S}(\mathbb{R}^{d})$ be the Schwartz class of rapidly decreasing function. Given $f \in \mathcal{S}(\mathbb{R}^{d})$, define the Fourier transform as
 \begin{equation*}
\hat{f}(\xi)=\mathcal{F}f(\xi)=(2\pi)^{-d/2}\int_{\mathbb{R}^{d}}e^{-ix\cdot\xi}f(x)\mathrm{d}x,
\end{equation*}
and its inverse Fourier transform as
\begin{equation*}
\check{f}(\xi)=\mathcal{F}^{-1}f(\xi)=(2\pi)^{-d/2}\int_{\mathbb{R}^{d}}e^{ix\cdot\xi}f(\xi)\mathrm{d}\xi.
\end{equation*}
Let $(\chi,\varphi)$ be a couple of smooth functions with values in $[0,1]$ such that $\chi$ is supported in the ball $\mathcal{B} :=\{\xi \in \mathbb{R}^{d}, |\xi| \leq \frac{4}{3} \}$, $\varphi$ is supported in the annulus $\mathcal{C} :=\{\xi \in \mathbb{R}^{d}, \frac{3}{4} \leq |\xi| \leq \frac{8}{3} \}$, $\varphi = \chi(\frac{\xi}{2})-\chi(\xi)$ and
\begin{align}
\chi(\xi)+\sum_{j \in \mathbb{N}} \varphi(2^{-j}\xi)=1,\quad \forall \xi \in \mathbb{R}^{d}; \quad \quad \sum_{j \in \mathbb{Z}}\varphi(2^{-j}\xi)=1,\quad \forall \xi \in \mathbb{R}^{d}\setminus \{0\}. \nn
\end{align}
For every $f \in \mathcal{S}'(\mathbb{R}^{d})$, we define the inhomogeneous Littlewood-Paley operators as follows
\begin{align}
\Delta_{-1}f: = \chi(D)u; \quad \Delta_{j}f: =\phi(2^{-j}D)f,\quad S_{j}f: = \sum_{-1 \leq k \leq j-1}\Delta_{k}u,\quad \forall j \in \mathbb{N}\cup\{0\}.\nn
\end{align}
The homogeneous Littlewood-Paley operators can be defined as follows
\begin{align}
\dot{\Delta}_{j}f: = \phi(2^{-j}D)f,\quad \dot{S}_{j}f: = \sum_{k \in \mathbb{Z},k\leq j-1}\Delta_{k}u, \quad \forall j \in \mathbb{Z}.\nn
\end{align}
Moreover, we have the following Bony's para-product decomposition
\begin{align}
uv = \sum_{j \geq -1}S_{j-1}u \Delta_{j}v + \sum_{j \geq -1}S_{j-1}v \Delta_{j}u + \sum_{j \geq -1} \Delta_{j}u \tilde{\Delta}_{j}v, \nn
\end{align}
where
\begin{align}
\tilde{\Delta}_{j}v = \Delta_{j-1}v + \Delta_{j}v + \Delta_{j+1}v. \nn
\end{align}
We now recall the definitions of the Besov type spaces.
\begin{definition}\label{def 2.1}
Let $s\in\mathbb{R}$ and $1\leq p,~q\leq+\infty$, the inhomogeneous Besov space $B_{p,q}^{s}(\mathbb{R}^{d})$ abbreviated as $B_{p,q}^{s}$ is defined by
\begin{equation*}
B_{p,q}^{s}=\{f(x)\in \mathcal{S}'(\mathbb{R}^{d});\|f\|_{B_{p,q}^{s}}<+\infty\},
\end{equation*}
where
\begin{equation*}
\|f\|_{B_{p,q}^{s}}=\begin{cases}
{ \begin{array}{ll} (\sum_{j\geq-1}2^{jsq}\|\Delta_{j} f\|_{L^{p}}^{q})^{\frac{1}{q}},~~~\mathrm{for}~q<+\infty,\\
\sup_{j\geq-1}2^{js}\|\Delta_{j} f\|_{L^{p}},~~~~~~~\mathrm{for}~q=+\infty
 \end{array} }
\end{cases}
\end{equation*}
is the inhomogeneous Besov norm.\\
The homogeneous Besov space $\dot{B}_{p,q}^{s}(\mathbb{R}^{d})$ abbreviated as $\dot{B}_{p,q}^{s}$ is defined by
\begin{equation*}
\dot{B}_{p,q}^{s}=\{f(x)\in \mathcal{S}'(\mathbb{R}^{d})/ \mathcal{P}(\mathbb{R}^{d});\|f\|_{\dot{B}_{p,q}^{s}}<+\infty\},
\end{equation*}
where $\mathcal{P}(\mathbb{R}^{d})$ represents the set of polynomial functions, and
\begin{equation*}
\|f\|_{\dot{B}_{p,q}^{s}}=\begin{cases}
{ \begin{array}{ll}(\sum_{j\in \mathbb{Z}}2^{jsq}\|\dot{\Delta}_{j} f\|_{L^{p}}^{q})^{\frac{1}{q}},~~~\mathrm{for}~1\leq q<+\infty,\\
\sup_{j\in \mathbb{Z}}2^{js}\|\dot{\Delta}_{j} f\|_{L^{p}},~~~~~~~\mathrm{for}~q=+\infty
 \end{array} }
\end{cases}
\end{equation*}
is the homogeneous Besov norm.\\
For $s \in \mathbb{R}$ and $1 \leq p,q \leq \infty$, the norms of Besov type spaces  $B_{p,q}^{s,g}$ and $\dot{B}_{p,q}^{s,g}$ are defined as follows
  \begin{align}
  \|f\|_{B_{p,q}^{s,g}}  = \|2^{js}g(2^j)\|\Delta_{j}f\|_{L^p}\|_{l^q} < \infty, \nn\\
   \|f\|_{\dot{B}_{p,q}^{s,g}}  = \|2^{js}g(2^j)\|\dot{\Delta}_{j}f\|_{L^p}\|_{l^q} < \infty, \nn
\end{align}
where $g$ is a non-negative non-decreasing smooth function.
\end{definition}
\par
The following lemma presents the Bernstein's inequalities.
\begin{lemma} \label{lem 2.1}
(\cite{QCZ2007,DLJD}) Let $\alpha \geq 0$. Let $1 \leq p \leq q \leq \infty$.\\
\hangindent=3em
\hangafter=1
 (1) If f satisfies
   \begin{align}
 supp \hat{f} \subset \{ \xi \in \mathbb{R}^{d}: |\xi| \leq K2^{j}\}\nn
\end{align}
for some integer j and a constant $K>0$, then
   \begin{align}
 \|(-\Delta)^{\alpha}f\|_{{L^q}(\mathbb{R}^{d})} \leq C_{1}2^{2 \alpha j +jd(\frac{1}{p}-\frac{1}{q})}\|f\|_{{L^p}(\mathbb{R}^{d})}. \nn
\end{align}
\hangindent=3em
\hangafter=1
 (2) If f satisfies
   \begin{align}
 supp \hat{f} \subset \{ \xi \in \mathbb{R}^{d}: K_{1}2^{j} \leq |\xi| \leq K_{2}2^{j}\} \nn
\end{align}
for some integer j and  constant $0 < K_{1} \leq K_{2}$, then
   \begin{align}
 C_{1}2^{2\alpha j}\|f\|_{{L^q}(\mathbb{R}^{d})} \leq \|(-\Delta)^{\alpha}f\|_{{L^q}(\mathbb{R}^{d})} \leq C_{2}2^{2 \alpha j +jd(\frac{1}{p}-\frac{1}{q})}\|f\|_{{L^p}(\mathbb{R}^{d})}, \nn
\end{align}
where $C_{1}$ and $C_{2}$ are constants depending on $\alpha$, p and q only.
\end{lemma}

In order to derive some lower bounds of the operator $\mathcal{L}$ which is defined by ({\ref{1.2}}), we can also write the operator $\mathcal{L}$ in a nonlocal dissipation operator form as:
\begin{equation}\label{L2}
\mathcal{L} f(x)=p.v.\int_{\mathbb{R}^2}\frac{f(x)-f(y)}{|x-y|^2g(1/|x-y|)|x-y|}dy,
\end{equation}
which is an equivalent definition, readers interested in the formula (\ref{L2}) may refer to \cite{MALV2014}. By using the equivalent form (\ref{L2}) of $\mathcal{L} f(x$, the following estimate can be easily deduced.
\begin{lemma} \label{lem 2.2}
(\cite{DLJD}) Let $\mathcal{L}$ be the operator defined by ({\ref{1.2}}) and (\ref{L2}). Then, for $1<p<\infty$,
    \begin{align}
      &|f(x)|^{p-2}f(x)(\mathcal{L}f(x)) \geq \frac{1}{p}\mathcal{L}(|f|^{p}). \nn
      \end{align}
\end{lemma}
By using Lemma 2.2, the following Lemma 2.3 can be proved easily, for details see \cite{DLJD}.
\begin{lemma}  \label{lem 2.3}
 Let $\mathcal{L}$ be the operator defined by ({\ref{1.2}}) and (\ref{L2}). Then, for $p \geq 2$,
    \begin{align}
      &\int_{\mathbb{R}^2} |f|^{p-2}f(\mathcal{L}f)\mathrm{d}x \geq \frac{2}{p}\int_{\mathbb{R}^2} |\mathcal{L}^{\frac{1}{2}}(|f|^{\frac{p}{2}})|^{2}\mathrm{d}x. \nn
      \end{align}
\end{lemma}
The following Lemma \ref{lem 2.4} is a generalized version of the Bernstein type inequality associated with the operator $\mathcal{L}$.

\begin{lemma} \label{lem 2.4}
(\cite{DLJD}) Let $j \geq 0$ be an integer  and $p \in [2,\infty)$. Let $\mathcal{L}$ be defined by ({\ref{1.2}}) and (\ref{L2}). Then, for any $f \in \mathcal{S}(\mathbb{R}^{2})$,
    \begin{equation}\label{2.1}
      2^{j}g^{-1}(2^j)\|\Delta_{j}f\|_{{L^p}(\mathbb{R}^{2})}^{p} \leq C \int_{\mathbb{R}^{2}}|\Delta_{j}f|^{p-2}\Delta_{j}f\mathcal{L}\Delta_{j}f\mathrm{d}x,
      \end{equation}
where C is a constant depending on p and d only.
\end{lemma}
\par
The following Lemma \ref{lem 2.5} will be used in the proof of many commutator estimates.
\begin{lemma} \label{lem 2.5}
(\cite{TSF2011}) Consider two different cases: $\delta \in (0, 1)$ and $\delta = 1$.\par
 (1) Let $\delta \in (0, 1)$ and $(p_{1},p_{2},p_{3})  \in [1,\infty]^{3}$. If $|x|^{\delta}h \in L^{1}(\mathbb{R}^2)$, $f \in \dot{B}_{p_{2},\infty}^{\delta}(\mathbb{R}^2)$ , $g \in L^{p_{3}}(\mathbb{R}^2)$ and $\frac{1}{p_{1}}=\frac{1}{p_{2}}+\frac{1}{p_{3}}$, then
     \begin{equation}\label{2.2}
      \|h \ast (fg)-f(h \ast g)\|_{L^{p_{1}}} \leq C\||x|^{\delta}h\|_{L^1}\|f\|_{\dot{B}_{p_{2},\infty}^{\delta}}\|g\|_{L^{p_{3}}},
      \end{equation} \par
where C is a constant independent of f, g and h.\par
 (2) Let $\delta = 1$ and $p \in [1,\infty]$. Let $r_{1} \in [1,p]$ and $r_{2} \in [1,\infty]$ satisfying $\frac{1}{r_{1}}+\frac{1}{r_{2}}=1$. Then
    \begin{equation}\label{2.3}
      \|h \ast (fg)-f(h \ast g)\|_{L^p} \leq C\||x|h\|_{L^{r_{1}}}\|\nabla f\|_{L^p}\|g\|_{L^{r_{2}}}.
      \end{equation}
\end{lemma} \par
Next we recall some classical commutator estimates.
\begin{lemma}\label{lem 2.6}
(\cite{KPV1991,KP1988}) Let $s>0$, $1<r<\infty$ and $\frac{1}{r}=\frac{1}{p_{1}}+\frac{1}{q_{1}}=\frac{1}{p_{2}}+\frac{1}{q_{2}}$ with $q_{1},p_{2}\in(1,\infty)$ and $p_{1},q_{2}\in[1,\infty]$. Then
\begin{align}
\|[\Lambda^{s},f]g\|_{L^{r}}\leq C(\|\nabla f\|_{L^{p_{1}}}\|\Lambda^{s-1}g\|_{L^{q_{1}}}+\|\Lambda^{s}f\|_{L^{p_{2}}}\|g\|_{L^{q_{2}}}),\nn
\end{align}
where $[\Lambda^{s},f]g=\Lambda^{s}(fg)-f\Lambda^{s}g$ and $C$ is a constant depending on the indices $s,r,p_{1},q_{1},p_{2},q_{2}$.
\end{lemma}

\begin{lemma} \label{lem 2.7}
 Let $\mathcal{R}_{g}$ be defined as $\mathcal{L}^{-1}\partial_{x_1}$. Assume
    \begin{align}
      &(p_{1},p_{2},p_{3}) \in [2,\infty)^{3}, \quad \frac{1}{p_1}=\frac{1}{p_2}+\frac{1}{p_3}, \quad q \in [1,\infty], \quad 0 < s < \delta. \nn
    \end{align}
Let  $[\mathcal{R}_{g},u]F = \mathcal{R}_{g}(uF)-u\mathcal{R}_{g}F$ be a standard commutator. Then
    \begin{align}
      &\|[\mathcal{R}_{g},u]F\|_{B_{p_{1},q}^{s,g}} \leq C(\|u\|_{\dot{B}_{p_{2},\infty}^{\delta}}\|F\|_{B_{p_{3},q}^{s-\delta,g^2}}+\|u\|_{L^2}\|F\|_{L^2}), \nn
      \end{align}
where C denotes a constant independent of $s$, $q$ and $p_{i}$ for $i=1,\ 2,\ 3$.
\end{lemma}
The Lemma \ref{lem 2.7} can be proven in a similar fashion as that of Proposition 2.6 in \cite{DLJD} by Lemma \ref{lem 2.5}, we omit the details for conciseness.\par
The following Lemma \ref{lem 2.8} is a general interpolation formula.
\begin{lemma} \label{lem 2.8}
(\cite{DLJD}) Let $\beta \in (2,\infty)$, $s \in (0,1)$, $0<\varepsilon(\beta-2) \leq 2$ and $f \in L^{\frac{2 \beta}{1+ \varepsilon}}(\mathbb{R}^2) \cap \dot{H}^{s+(1-\frac{2}{\beta})(1+\varepsilon)}(\mathbb{R}^2)$.
 Then
   \begin{equation}\label{2.4}
 \||f|^{\beta -2}f\|_{\dot{H}^{s}} \leq C\|f\|_{L^{\frac{2\beta}{1+\varepsilon}}}^{\beta-2}\|f\|_{\dot{B}_{\frac{2\beta}{2-\varepsilon(\beta-1)},2}^{s}} \leq C\|f\|_{L^{\frac{2\beta}{1+\varepsilon}}}^{\beta-2}\|f\|_{\dot{H}^{s+(1-\frac{2}{\beta})(1+\varepsilon)}}.
      \end{equation}
\end{lemma}
\begin{lemma} \label{lem 2.9}
(\cite{DJ2012}) Let $T > 0$ and u be a divergence-free smooth vector field satisfying
\begin{align}
\int_{0}^{T} \|\nabla u \|_{L^\infty}\mathrm{d}t < \infty. \nn
\end{align}
Assume that $\theta$ solves
\begin{align}
\partial_{t}\theta + u \cdot \nabla \theta =f. \nn
\end{align}
Let $g : (0, \infty) \rightarrow (0, \infty)$ be an nondecreasing and radially symmetric function satisfying ({\ref{1.3}}). Let $\theta_0(x)\in B^{0,g}_{p,1}(\mathbb{R}^2)$, $f(t,x)\in L^1([0,T]; B^{0,g}_{p,1}(\mathbb{R}^2))$ and  $p \in [1,\infty]$. Then for any $t > 0$,
\begin{align}
\|\theta\|_{B_{p,1}^{0,g}} \leq (\|\theta_{0}\|_{B_{p,1}^{0,g}}+\|f\|_{L^{1}_{t}{B_{p,1}^{0,g}}})(1+\int_{0}^{t} \|\nabla u\|_{L^\infty}\mathrm{d}t).\nn
\end{align}
\end{lemma}
\setcounter{equation}{0}

    \section{some $\textit{A priori}$ estimations on $u$ and $\theta$}
    \setcounter{equation}{0}
\vskip .1in
\quad This section is divided into four subsections. The first subsection mainly proves $\|u\|_{\dot{H}^{1}}$ and $\|\theta\|_{L^\infty}$. In the second subsection, we prove the estimate of the $\|G\|_{L^r}$ and $\|G\|_{L^{1}_{t}{B_{r,1}^{\frac{2}{r}}}}$ for $r\in(2,4)$. In the third subsection, we further derive the estimates of  $\|\omega\|_{L^{1}_{t}{B_{\infty,1}^{0}}}$, $\|\nabla u\|_{L_{t}^{1}L^\infty}$ and $\|\nabla \theta\|_{L^\infty}$. The final subsection establishes the global $H^{s}$ estimates of $u$ and $\theta$ for $s\ge 2$.
\subsection{The estimates of $\|G\|_{L^{2}}$ and $\|\theta\|_{L^{\infty}}$}
     \quad In this subsection, we first prove the global estimates of $\|u\|_{L^2}$ and $\|\theta\|_{L^p}$ with $p \in [2,4)$, then by using the commutator estimation of Lemma {\ref{lem 2.7}}, Young's and Gr$\rm\ddot{o}$nwall's inequalities and interpolation inequality, we get the estimate of $\|G\|_{L^2}$, and  further show the estimate of $\|\theta\|_{L^\infty}$.
\begin{proposition} \label{pro 3.1}
Assume $(u_{0},\theta_{0})$ satisfies the assumptions in Theorem ${\ref{thm 1.1}}$.  Then the strong solution $(u,\theta)$ of ({\ref{1.1}}) admits the following estimates, for $0<t<T$,
 \begin{equation} \label{3.1}
    \|u(t)\|_{L^2}^{2}+\int_{0}^{t}\|\mathcal{L}^\frac{1}{2}u(\tau)\|_{L^2}^{2} \mathrm{d}\tau \leq (\|u_{0}\|_{L^2}^{2}+\|\theta_{0}\|_{L^2}^{2})e^{CT},
    \end{equation}
     \begin{equation} \label{3.2}
    \|\theta\|_{L^p} \leq (\|u_{0}\|_{L^2}^{2}+\|\theta_{0}\|_{L^2}^{2})e^{CT} , \forall\ p \in [2,4) .
    \end{equation}
\end{proposition}
\begin{proof}
Taking the $L^2$ inner product of the velocity equation with $u$ and the temperature equation with $\theta$, we find by adding them up
    \begin{align}
    &\frac{1}{2}\frac{\mathrm{d}}{\mathrm{d}t}(\|u\|_{L^2}^{2}+\|\theta\|_{L^2}^{2})+\|\mathcal{L}^\frac{1}{2}u\|_{L^2}^{2}
    \leq 2(\|u\|_{L^2}^{2}+\|\theta\|_{L^2}^{2}),\nonumber
    \end{align}
which together with Gr$\rm\ddot{o}$nwall inequality implies that
    \begin{equation}\label{3.3}
    \|u\|_{L^2}^{2}+\|\theta\|_{L^2}^{2}+\int_{0}^{t}\|\mathcal{L}^\frac{1}{2}u(\tau)\|_{L^2}^{2} \mathrm{d}\tau  \leq (\|u_{0}\|_{L^2}^{2}+\|\theta_{0}\|_{L^2}^{2})e^{CT}.
    \end{equation}
Then, multiplying the temperature equation by $|\theta|^{p-2}\theta$ and integrating the resulting equality with respect to $x$, we get
  \begin{align}
    &\frac{1}{p}\frac{\mathrm{d}}{\mathrm{d}t}\|\theta\|_{L^p}^{p}
    =\int_{\mathbb{R}^{2}}u_{2}|\theta|^{p-2}\theta \mathrm{d}x
     \leq \|u\|_{L^p}\|\theta\|_{L^p}^{p-1}.\nn
    \end{align}
By ({\ref{3.3}}) and a interpolation of $L^p$ between $L^2$ and $\dot{H}^{\frac{1}{2}-}$, we have
  \begin{align}
    \|\theta\|_{L^p}   \nonumber
    &\leq \|\theta_{0}\|_{L^p}+\int_{0}^{t}\|u\|_{L^p}\mathrm{d}\tau\nn\\
     &\leq  \|\theta_{0}\|_{L^p}+\int_{0}^{t}(\|u\|_{L^2}+\|\mathcal{L}^{\frac{1}{2}}u\|_{L^2})\mathrm{d}\tau \nn\\
     &\leq (\|u_{0}\|_{L^2}^{2}+\|\theta_{0}\|_{L^2}^{2})e^{CT}.\nn
    \end{align}
The proof of Proposition {\ref{pro 3.1}} is complete.
\end{proof} \par
Next, we establish the following estimates of $\|G\|_{L^2}$ and $\|\theta\|_{L^\infty}$.

\begin{proposition} \label{pro 3.2}
Assume $(u_{0},\theta_{0})$ satisfies the assumptions in Theorem {\ref{thm 1.1}}.  Let $(u,\theta)$ be the strong solution  of ({\ref{1.1}}), then the quantity G and the temperature $\theta$ have the following estimates for any $t \in (0,T)$,
 \begin{equation} \label{3.4}
    \|G(t)\|_{L^2}^{2}+\int_{0}^{t}\|\mathcal{L}^\frac{1}{2}G(\tau)\|_{L^2}^{2} \mathrm{d}\tau \leq C(T,u_{0},\theta_{0}) ,
    \end{equation}
     \begin{equation}\label{3.x}
   \int_{0}^{t}\|\omega\|_{B_{p,\infty}^{0,g^{-1}}}^{2} \mathrm{d}\tau \leq C(T,u_{0},\theta_{0}),\ 2\le p<4,
    \end{equation}
     \begin{equation} \label{3.5}
    \|\theta\|_{L^\infty} \leq C(T,u_{0},\theta_{0}),
    \end{equation}
 where and in the sequel $C(T,u_{0},\theta_{0})$ denotes a constant depending the $H^s-$norm of the initial data $(u_0,\theta_0)$ and $T$.
\end{proposition}

\begin{proof}
Multiplying ({\ref{1.5}}) by G, using the integration by parts, we obtain
    \begin{equation}\label{3.6}
    \frac{1}{2}\frac{\mathrm{d}}{\mathrm{d}t}\|G\|_{L^2}^{2}+\|\mathcal{L}^\frac{1}{2}G\|_{L^2}^{2}
    = \int_{\mathbb{R}^{2}}\mathcal{R}_{g}u_{2}G\mathrm{d}x+\int_{\mathbb{R}^{2}}[\mathcal{R}_{g},u \cdot \nabla]\theta G \mathrm{d}x.
    \end{equation}
It follows from the H$\rm\ddot{o}$lder inequality and the interpolation inequality that
    \begin{equation}\label{3.7}
    \begin{split}
    \lvert \int_{\mathbb{R}^{2}} \mathcal{R}_{g}u_{2}G\mathrm{d}x \rvert
     &\leq \|u\|_{L^2}\|\mathcal{R}_{g}G\|_{L^2} \\
     &\leq C\|u\|_{L^2}(\|G\|_{L^2}+\|\mathcal{L}^\frac{1}{2}G\|_{L^2}) \\
     &\leq C\|u\|_{L^2}^{2}+C\|G\|_{L^2}^{2}+\frac{1}{2}\|\mathcal{L}^\frac{1}{2}G\|_{L^2}^{2}.
     \end{split}
    \end{equation}
Next by  the H$\rm\ddot{o}$lder inequality, we take $\varepsilon \in (0,\frac{1}{2})$ small enough,
 \begin{align}
    \lvert \int_{\mathbb{R}^{2}}[\mathcal{R}_{g},u \cdot \nabla]\theta G \mathrm{d}x \rvert
    =\lvert \int_{\mathbb{R}^{2}}G  \nabla \cdot [\mathcal{R}_{g},u]\theta  \mathrm{d}x \rvert
    \leq C\|G\|_{\dot{H}^{\frac{1}{2}-\varepsilon}}\|[\mathcal{R}_{g},u] \theta \|_{\dot{H}^{\frac{1}{2}+\varepsilon}}.\nn
    \end{align}
In addition, due to the condition in ({\ref{1.3}}),
\begin{equation}\label{3.8}
    \|G\|_{\dot{H}^{\frac{1}{2}-\varepsilon}}^{2} \leq \|G\|_{H^{\frac{1}{2}-\varepsilon}}^{2} = \sum_{j \geq -1} 2^{j-2 \varepsilon j}\|\Delta_{j}G\|_{L^2}^{2} \leq C\sum_{j \geq -1} 2^{j}g^{-2}(2^j)\|\Delta_{j}G\|_{L^2}^{2} \leq C\|\mathcal{L}^{\frac{1}{2}}G\|_{L^2}^{2}.
    \end{equation}
Applying Lemma {\ref{lem 2.7}} with $  \frac{1}{2}<\delta=\frac{1}{2}+\varepsilon<1 $ and $2<p<4$, we obtain
\begin{align}
    \|[\mathcal{R}_{g},u] \theta \|_{\dot{H}^{\frac{1}{2}+\varepsilon}}
      \leq \|[\mathcal{R}_{g},u] \theta \|_{B^{\frac{1}{2}+\varepsilon}_{2,2}}
      \leq C\|u\|_{B_{\frac{2p}{p-2},\infty}^{2\varepsilon+\frac{1}{2}}} \|\theta\|_{B_{p,2}^{-{\varepsilon},g}}+\|u\|_{L^2} \|\theta\|_{L^2}.\nn
    \end{align}
Since $u={\nabla}^{\perp} {\Delta}^{-1} \omega$,
   \begin{align}
    \|u\|_{B_{\frac{2p}{p-2},\infty}^{2\varepsilon+\frac12}}
     &\leq \sup_{j \geq -1} 2^{j(2\varepsilon+\frac12)}\|\Delta_{j}u\|_{L^{\frac{2p}{p-2}}} \nn\\
     &\leq \|\Delta_{-1}u\|_{L^{\frac{2p}{p-2}}}+\sup_{j \geq 0} 2^{j(2\varepsilon-\frac12)}\|\Delta_{j}\omega\|_{L^{\frac{2p}{p-2}}} \nn  \\
     &\leq \|u\|_{L^2}+\sup_{j \geq 0} 2^{j(2\varepsilon-\frac12)}\|\Delta_{j}G\|_{L^{\frac{2p}{p-2}}}+\sup_{j \geq 0} 2^{j(2\varepsilon-\frac12)}\|\Delta_{j}\mathcal{R}_{g}\theta\|_{L^{\frac{2p}{p-2}}} \nn \\
     &\leq \|u\|_{L^2}+\sup_{j \geq 0} 2^{j(2\varepsilon-\frac12+\frac{2}{p})}\|\Delta_{j}G\|_{L^2}+\sup_{j \geq 0} 2^{j(2\varepsilon-\frac32+\frac{4}{p})}g(2^j)\|\Delta_{j}\theta\|_{L^p}. \nn
      \end{align}
Choosing $p \in (2,4)$ such that  $s:=2\varepsilon-\frac12+\frac{2}{p}<\frac12$ and $2\varepsilon-\frac32+\frac{4}{p}<0$, for example $p=3$, we have
\begin{align}
&\sup_{j \geq 0} 2^{j(2\varepsilon-\frac12+\frac{2}{p})}\|\Delta_{j}G\|_{L^2} \leq \|G\|_{H^s} \leq C\|G\|_{L^2}^{\lambda}\|\mathcal{L}^{\frac{1}{2}}G\|_{L^2}^{1-\lambda},\nn
\end{align}
where $\lambda \in (0,1)$. Clearly $\|\theta\|_{B_{p,2}^{-\varepsilon,g}} \leq \|\theta\|_{L^p}$, therefore
   \begin{equation}\label{3.9}
    \|[\mathcal{R}_{g},u] \theta \|_{\dot{H}^{\frac{1}{2}+\varepsilon}} \leq C(\|u\|_{L^2}+\|G\|_{L^2}^{\lambda}\|\mathcal{L}^{\frac{1}{2}}G\|_{L^2}^{1-\lambda}+\|\theta\|_{L^p})\|\theta\|_{L^p}+\|u\|_{L^2}\|\theta\|_{L^2}.
    \end{equation}
Substituting (\ref{3.7}), (\ref{3.8}) and (\ref{3.9}) into ({\ref{3.6}}), we obtain
   \begin{equation}
   \begin{split}
   \frac{\mathrm{d}}{\mathrm{d}t} \|G\|_{L^2}^{2} + 2\|\mathcal{L}^{\frac{1}{2}}G\|_{L^2}^{2} &\leq C\|u\|_{L^2}\|\theta\|_{L^p}\|\mathcal{L}^{\frac{1}{2}}G\|_{L^2}+C\|\theta\|_{L^p}\|G\|_{L^2}^{\lambda}\|\mathcal{L}^{\frac{1}{2}}G\|_{L^2}^{2-\lambda}\\&+C\|\theta\|_{L^p}^{2}\|\mathcal{L}^{\frac{1}{2}}G\|_{L^2}+C\|u\|_{L^2}\|\theta\|_{L^2}\|\mathcal{L}^{\frac{1}{2}}G\|_{L^2}\\
   &+C\|u\|_{L^2}^{2}+C\|G\|_{L^2}^{2}+\frac{1}{2}\|\mathcal{L}^\frac{1}{2}G\|_{L^2}^{2}, \nn
   \end{split}
    \end{equation}
which together with the Young's inequality and Gr$\rm\ddot{o}$nwall's inequality implies the desired estimate ({\ref{3.4}}). Recalling that $G=\omega -\mathcal{R}_{g}\theta$, we have for $2\le p<4$
   \begin{align}
   \int_{0}^{t}\|\omega\|_{B_{p,\infty}^{0,g^{-1}}}^{2} \mathrm{d}\tau &\leq \int_{0}^{t}\|G\|_{B_{p,\infty}^{0,g^{-1}}}^{2} \mathrm{d}\tau + \int_{0}^{t}\|\mathcal{R}_{g}\theta\|_{B_{p,\infty}^{0,g^{-1}}}^{2} \mathrm{d}\tau \nn\\
   &\leq \int_{0}^{t}\|G\|_{L^p}^{2} \mathrm{d}\tau + \int_{0}^{t}\|\theta\|_{L^p}^{2} \mathrm{d}\tau \nn \\
   &\leq C(T,u_{0},\theta_{0}),\nn
    \end{align}
thus
    \begin{align}
     \int_{0}^{t} \|u(\tau)\|_{L^\infty}\mathrm{d}\tau &\leq \int_{0}^{t}\|u\|_{L^2}^{\lambda} \|u\|_{W^{1-\varepsilon,p}}^{1-\lambda} \mathrm{d}\tau \nn\\
     &\leq \int_{0}^{t}\|u\|_{L^2}^{\lambda} (\sum_{j \geq -1} 2^{j(1-\varepsilon)}\|\Delta_{j}u\|_{L^p})^{1-\lambda} \mathrm{d}\tau \nn\\
          &\leq \int_{0}^{t}\|u\|_{L^2}^{\lambda} (\|u\|_{L^2}+\sum_{j \geq 0} 2^{-\varepsilon j}g(2^j)g^{-1}(2^j)\|\Delta_{j}\omega\|_{L^p})^{1-\lambda} \mathrm{d}\tau \nn\\
      &\leq C\int_{0}^{t}\|u\|_{L^2}^{\lambda}(1+ \|\omega\|_{B_{p,\infty}^{0,g^{-1}}})^{1-\lambda}\mathrm{d}\tau \nn\\
     &\leq C(T,u_{0},\theta_{0}) ,\nonumber
    \end{align}
where $\lambda \in (0,1)$. We thus get from the temperature equation
\begin{align}
    &\|\theta\|_{L^\infty} \leq C\|\theta_{0}\|_{L^\infty}+\int_{0}^{t} \|u(\tau)\|_{L^\infty}\mathrm{d}\tau \leq C(T, u_{0}, \theta_{0}), \nonumber
    \end{align}
and  the proof of Proposition {\ref{pro 3.2}} is completed.
\end{proof}

\subsection{The estimates of $\|G\|_{L^{r}}$ and $\|G\|_{L^{1}_{t}{B_{r,1}^{\frac{2}{r}}}}$}
 In this subsection, we first establish the estimate of $\|G\|_{L^r}$ for $r \in (2,4)$, and further prove the estimate of $\|G\|_{L^{1}_{t}{B_{r,1}^{\frac{2}{r}}}}$, which serves as an important step towards a global bound for $\|\nabla u\|_{L_{t}^{1}{L^\infty}}$.
\begin{proposition} \label{pro 3.3}
Assume $(u_{0},\theta_{0})$ satisfies the assumptions in Theorem {\ref{1.1}}. Let $(u, \theta)$ be a strong solution of (\ref{1.1}), then, for any $r\in(2,4)$, $G$ obeys the following the bound, for any $0<t<T$,
    \begin{align}
    \|G\|_{L^r}^{r}+ C\int_{0}^{t} \|G\|_{L^{\frac{2r}{1+\varepsilon}}}^{r} \mathrm{d} \tau \leq C(T, u_{0}, \theta_{0}). \nonumber
    \end{align}
\end{proposition}
\begin{proof}
Multiplying $({\ref{1.5}})$ by $G|G|^{r-2}$ and integrating with respect to $x$, we obtain
    \begin{equation}\label{3.10}
    \frac{1}{r} \frac{\mathrm{d}}{\mathrm{d}t}\|G\|_{L^r}^{r}+\int_{\mathbb{R}^{2}}G|G|^{r-2}\mathcal{L}G\mathrm{d}x
    =\int_{\mathbb{R}^{2}}[\mathcal{R}_{g},u \cdot \nabla]\theta|G|^{r-2}G\mathrm{d}x-\int_{\mathbb{R}^{2}} \mathcal{R}_{g}u_{2}|G|^{r-2}G\mathrm{d}x.
    \end{equation}
By Lemma {\ref{lem 2.4}}, we have
\begin{align}
    &\int_{\mathbb{R}^{2}}(\mathcal{L}G)|G|^{r-2}G\mathrm{d}x \geq C \int \lvert \mathcal{L}^{\frac{1}{2}}(|G|^{\frac{r}{2}}) \rvert^{2} \mathrm{d}x=C\|\mathcal{L}^{\frac{1}{2}}(|G|^{\frac{r}{2}})\|_{L^2}^{2}. \nn
    \end{align}
Choosing $\varepsilon > 0$ small enough such that
\begin{align}
   (1+\varepsilon)(1- \frac{2}{r}) < \frac{1}{2}. \nn
    \end{align}
Thanks to the condition in ({\ref{1.3}}) and by a Sobolev embedding $\dot{H}^{\frac{1-\varepsilon}{2}}\hookrightarrow L^{\frac{4}{\varepsilon+1}}$, we obtain
\begin{equation}\label{3.11}
\begin{split}
   \|\mathcal{L}^{\frac{1}{2}}(|G|^{\frac{r}{2}})\|_{L^2}^{2}  &= \sum_{j \geq -1} \|\Delta_{j} \mathcal{L}^{\frac{1}{2}}(|G|^{\frac{r}{2}})\|_{L^2}^{2} \\
   &\geq \sum_{j \geq -1} 2^{j} g^{-1}(2^j)\|\Delta_{j} (|G|^{\frac{r}{2}})\|_{L^2}^{2}  \\
   &\geq c\sum_{j \geq -1} 2^{(1-\varepsilon)j}\|\Delta_{j} (|G|^{\frac{r}{2}})\|_{L^2}^{2}  \\
   &=c\|\Lambda^{{\frac{1}{2}}-{\frac{\varepsilon}{2}}}(|G|^{\frac{r}{2}})\|_{L^2}^{2}  \\
   &\geq c \|G\|_{L^{\frac{2r}{1+\varepsilon}}}^{r},
   \end{split}
    \end{equation}
where $c$ is a positive constant. For $r \in (2,4)$, we choose $s > 0$ small enough such that
\begin{align}
  s > \varepsilon,\quad s+(1+\varepsilon)(1- \frac{2}{r}) =\frac{1}{2}-\varepsilon. \nn
    \end{align}
We estimate the first term in the right hand side of ({\ref{3.10}}) as follows
\begin{align}
  \lvert \int_{\mathbb{R}^{2}}  \nabla \cdot [\mathcal{R}_{g},u]\theta G|G|^{r-2}  \mathrm{d}x \rvert
    \leq \|G|G|^{r-2}\|_{\dot{H}^{s}}\|[\mathcal{R}_{g},u] \theta \|_{\dot{H}^{1-s}}.\nn
    \end{align}
By Lemma {\ref{lem 2.8}}, we have
\begin{equation}\label{3.12}
\begin{split}
 \|G|G|^{r-2}\|_{\dot{H}^{s}} &\leq C\|G\|_{L^{\frac{2r}{1+\varepsilon}}}^{r-2}\|G\|_{\dot{H}^{s+(1+\varepsilon)(1- \frac{2}{r})}} \\&= C\|G\|_{L^{\frac{2r}{1+\varepsilon}}}^{r-2}\|G\|_{\dot{H}^{\frac{1}{2}-\varepsilon}} \\ &\leq C\|G\|_{L^{\frac{2r}{1+\varepsilon}}}^{r-2}\|\mathcal{L}^{\frac{1}{2}}G\|_{L^2}\\ &\leq \frac{c}{2}\|G\|^r_{L^{\frac{2r}{1+\varepsilon}}}+C\|\mathcal{L}^{\frac{1}{2}}G\|_{L^2}^{\frac{r}{2}}.
 \end{split}
    \end{equation}
By Lemma {\ref{lem 2.7}}, recalling $s > \varepsilon$ and $u={\nabla}^{\perp} {\Delta}^{-1} \omega$
\begin{equation}\label{3.13}
\begin{split}
  \|[\mathcal{R}_{g},u] \theta \|_{\dot{H}^{1-s}} &\leq C \|u\|_{\dot{B}_{2,\infty}^{1-s+\varepsilon}}\|\theta\|_{B_{\infty,2}^{-\varepsilon,g}}+C\|u\|_{L^2}\|\theta\|_{L^2} \\
  &\leq C\|\omega\|_{B_{2,2}^{0, g^{-1}}}\|\theta\|_{L^{\infty}}+C\|u\|_{L^2}\|\theta\|_{L^2}\\
  &\leq C(\|G\|_{L^2}+\|\theta\|_{L^2})\|\theta\|_{L^\infty}+C\|u\|_{L^2}\|\theta\|_{L^2}\\
  &\leq C(T,u_{0},\theta_{0}).
  \end{split}
    \end{equation}
The second term in the right hand side of ({\ref{3.10}}) admits the following estimate
\begin{equation}\label{3.14}
\begin{split}
  \lvert \int_{\mathbb{R}^{2}}G|G|^{r-2}  \mathcal{R}_{g}u_{2}  \mathrm{d}x \rvert
   & \leq C \|G\|_{L^r}^{r-1}\|\mathcal{R}_{g}u_{2} \|_{L^r}\\
   & \leq C\|G\|_{L^r}^{r}+C\|\mathcal{R}_{g}u_{2}\|_{L^r}^{r}.
   \end{split}
    \end{equation}
In view of the Littlewood-Paley decomposition and the Bernstein inequality (see Lemma {\ref{lem 2.1}}), we obtain
 \begin{equation}\label{3.15}
 \begin{split}
 \|\mathcal{R}_{g}u_{2}\|_{L^r} &\leq \sum_{j \geq -1}g(2^j)\|\Delta_{j}u\|_{L^r} \\
 &\leq \sum_{j \geq -1} g(2^j) \|\Delta_{j}u\|_{L^2}^{\frac{2}{r}}\|\Delta_{j}u\|_{L^\infty}^{1-\frac{2}{r}}\\
 &\leq C\sum_{j \geq -1}g(2^j)(2^{j\frac{2}{r}}\|\Delta_{j}u\|_{L^r})^{1-\frac{2}{r}}\\
 &\leq C+C\sum_{j \geq0}g(2^j)2^{-j(1-\frac{2}{r})^{2}j}\|\Delta_{j}\omega\|_{L^r}^{1-\frac{2}{r}}\\
  &\leq C+C\sum_{j \geq0}2^{-j\frac{1}{2}(1-\frac{2}{r})^{2}}\|\Delta_{j}\omega\|_{L^r}^{1-\frac{2}{r}}\\
  &\leq C+C\|\omega\|_{B_{r,\infty}^{-\varepsilon}}^{1-\frac{2}{r}},
 \end{split}
 \end{equation}
where $\varepsilon=\frac{1}{4}(1-\frac{2}{r})<1$, thus we have
  \begin{align}
  \int_{0}^{t}\|\mathcal{R}_{g}u_{2}\|_{L^r}^{r}\mathrm{d}\tau \leq \int_{0}^{t}C\|\omega\|_{B_{r,\infty}^{-\varepsilon}}^{r-2}\mathrm{d}\tau \leq C(T,u_{0},\theta_{0}).\nn
  \end{align}
 Putting (\ref{3.11}), (\ref{3.12}), (\ref{3.13}), (\ref{3.14}) and (\ref{3.15}) into (\ref{3.10}), by Young's inequality, we obtain
    \begin{align}
    \frac{1}{r} \frac{\mathrm{d}}{\mathrm{d}t}\|G\|_{L^r}^{r} +C\|G\|_{L^{\frac{2r}{1+\varepsilon}}}^{r}  &\leq C\|G\|_{L^r}^{r}+C\|\mathcal{L}^{\frac{1}{2}}G\|_{L^2}^{\frac{r}{2}} +C\|\omega\|_{B_{r,\infty}^{-\varepsilon}}^{r-2}+C. \nonumber
    \end{align}
Applying Gr$\rm\ddot{o}$nwall's inequality to the above inequality, we have
    \begin{align}
    &\|G\|_{L^r}^{r} +C\int_{0}^{t}\|G\|_{L^{\frac{2r}{1+\varepsilon}}}^{r}\mathrm{d}\tau \leq C(T,u_{0},\theta_{0}),\nn
    \end{align}
which completes the proof of Proposition {\ref{pro 3.3}}.
\end{proof}

\begin{proposition} \label{pro 3.4}
Assume that the initial data $(u_{0},\theta_{0})$ satisfies the assumptions in Theorem {\ref{thm 1.1}}.  Then the strong solution $(u,\theta)$ of ({\ref{1.1}}) obeys the following bound, for any $0<t<T$ and $2<r<4$,
    \begin{equation}\label{3.16}
    \|G\|_{{L^1_t}{B_{r,1}^{\frac{2}{r}}}} \leq C(T, u_{0}, \theta_{0}).
    \end{equation}
\end{proposition}
\begin{proof}
Applying  the operator $\Delta_{j}$  to ({\ref{1.5}}) for $j \geq -1$, we obtain
  \begin{equation}\label{3.17}
   \partial_t \Delta_{j}G+(u \cdot \nabla)\Delta_{j}G + \mathcal{L}\Delta_{j}G
   =\Delta_{j}[\mathcal{R}_{g},u \cdot \nabla]\theta - [\Delta_{j},u \cdot \nabla]G - \Delta_{j}\mathcal{R}_{g}u_{2}.
    \end{equation}
 For notational convenience, we denote
\begin{align}
   &f_{j}:=\Delta_{j}[\mathcal{R}_{g},u \cdot \nabla]\theta-[\Delta_{j},u \cdot \nabla]G-\Delta_{j}\mathcal{R}_{g}u_{2}.\nonumber
    \end{align}
Multiplying ({\ref{3.17}}) by $|\Delta_{j}G|^{r-2}\Delta_{j}G$, using the divergence-free condition of $u$, we derive
\begin{equation}\label{3.18}
   \frac{1}{r}\frac{\mathrm{d}}{\mathrm{d}t}\|\Delta_{j}G\|_{L^r}^{r}+\int_{\mathbb{R}^{2}}(\mathcal{L}\Delta_{j}G)|\Delta_{j}G|^{r-2}\Delta_{j}G\mathrm{d}x=\int_{\mathbb{R}^{2}}f_{j}|\Delta_{j}G|^{r-2}\Delta_{j}G\mathrm{d}x.
    \end{equation}
According to Lemma {\ref{lem 2.4}} , for $j \geq 0$, the dissipation part can be bounded below by
\begin{equation}\label{3.19}
   \int_{\mathbb{R}^{2}}(\mathcal{L}\Delta_{j}G)|\Delta_{j}G|^{r-2}\Delta_{j}G\mathrm{d}x \geq C2^j g^{-1}(2^j)\|\Delta_{j}G\|_{L^r}^{r}.
    \end{equation}
Inserting ({\ref{3.19}}) into ({\ref{3.18}}), we obtain
\begin{align}
   &\frac{1}{r}\frac{\mathrm{d}}{\mathrm{d}t}\|\Delta_{j}G\|_{L^r}^{r}+C2^{j}g^{-1}(2^j)\|\Delta_{j}G\|_{L^r}^{r}
   \leq C\|f_j\|_{L^r}\|\Delta_{j}G\|_{L^r}^{r-1}.\nn
    \end{align}
 Due to ({\ref{1.3}}), taking $\varepsilon$ small enough, integrating the above inequality with respect to time, it yields
\begin{equation}\label{3.20}
   \|\Delta_{j}G\|_{L^r}
   \leq Ce^{-C2^{(1-\varepsilon)j}t}\|\Delta_{j}G(0)\|_{L^r}+C\int_{0}^{t}e^{-C2^{(1-\varepsilon)j}(t-\tau)}\|f_j\|_{L^r}\mathrm{d}\tau.
    \end{equation}
 Integrating ({\ref{3.20}}) over $[0,t]$ and using the Young's inequality, it yields
\begin{equation}\label{3.21}
\begin{split}
   \int_{0}^{t}\|\Delta_jG\|_{L^r}\mathrm{d}\tau
   &\leq C2^{-(1-\varepsilon)j}\|\Delta_{j}G(0)\|_{L^r}+C2^{-(1-\varepsilon)j}\int_{0}^{t}\|f_j\|_{L^r}\mathrm{d}\tau \\ &\leq C2^{-(1-\varepsilon)j}\Big\{\|\Delta_{j}G(0)\|_{L^r}+\int_{0}^{t}\|\Delta_{j}[\mathcal{R}_{g},u \cdot \nabla]\theta\|_{L^r}\mathrm{d}\tau +\int_{0}^{t}\|[\Delta_{j},u \cdot \nabla]G\|_{L^r}\mathrm{d}\tau \\ &+ \int_{0}^{t}\|\Delta_{j}\mathcal{R}_{g}u_{2}\|_{L^r}\mathrm{d}\tau\Big\}.
     \end{split}
       \end{equation}
On the other hand, for the low frequency part, we immediately obtain
\begin{align}
\int_{0}^{t} \|\Delta_{-1}G\|_{L^r}\mathrm{d} \tau \leq C\int_{0}^{t} \|G\|_{L^2}\mathrm{d} \tau \leq C(T).\nn
\end{align}
Combining the upper high-low frequency estimations, multiplying them by $2^{j\frac{2}{r}}$ and summing  $j$ from -1 to $\infty$, we obtain
\begin{equation}\label{3.22}
\begin{split}
   \int_{0}^{t}\|G\|_{B_{r,1}^{\frac{2}{r}}}\mathrm{d}\tau
   &\leq C\|G(0)\|_{B_{r,1}^{\frac{2}{r}+\varepsilon-1}}+C\int_{0}^{t}\|\mathcal{R}_{g}u_{2}\|_{B_{r,1}^{\frac{2}{r}+\varepsilon-1}}\mathrm{d}\tau+C\int_{0}^{t}\|[\mathcal{R}_{g},u \cdot \nabla]\theta(\tau)\|_{B_{r,1}^{\frac{2}{r}+\varepsilon-1}}\mathrm{d}\tau \\ &+C\int_{0}^{t} \sum_{j=-1}^{\infty}2^{j(\frac{2}{r}+\varepsilon-1)}\|[\Delta_{j},u \cdot \nabla]G]\|_{L^r}\mathrm{d}\tau+C(T).
   \end{split}
       \end{equation}
 We estimate the second term in the right hand side of ({\ref{3.22}}) as follows
  \begin{equation}\label{3.23}
  \begin{split}
   C\int_{0}^{t}\|\mathcal{R}_{g}u_{2}\|_{B_{r,1}^{\frac{2}{r}+\varepsilon-1}}\mathrm{d}\tau  &= C\int_{0}^{t}\sum \limits_{j \geq -1}2^{j(\frac{2}{r}+\varepsilon-1)}\|\Delta_{j}\mathcal{R}_{g}u_{2}\|_{L^r}\mathrm{d}\tau  \\
   &\leq C\int_{0}^{t}(C\|\Delta_{-1}u\|_{L^r}+\sum \limits_{j \geq 0}2^{j(\frac{2}{r}+\varepsilon-2)}g(2^j)\|\Delta_{j}\omega\|_{L^r})\mathrm{d}\tau  \\
   &\leq C\int_{0}^{t}(\|u\|_{L^2}+\|\omega\|_{B_{r,1}^{\frac{2}{r}+2\varepsilon-2}})\mathrm{d}\tau  \\
   &\leq C\int_{0}^{t}(\|u\|_{L^2}+\|G\|_{L^r}+\|\theta\|_{L^r})\mathrm{d}\tau  \\
   &\leq C(T,u_{0},\theta_{0}).
 \end{split}
 \end{equation}
For the third term in the right hand side of ({\ref{3.22}}), by Lemma {\ref{lem 2.7}} and choosing $\varepsilon$ small enough such that $\frac{2}{r}+2\varepsilon <1$, we have
  \begin{equation}\label{3.24}
  \begin{split}
   C\int_{0}^{t}\|[\mathcal{R}_{g},u \cdot \nabla]\theta\|_{B_{r,1}^{\frac{2}{r}+\varepsilon-1}} \mathrm{d}\tau
   &\leq C\int_{0}^{t}\|[\mathcal{R}_{g},u]\theta\|_{B_{r,1}^{\frac{2}{r}+\varepsilon}} \mathrm{d}\tau\\
   &\leq C\int_{0}^{t}\Big(\|u\|_{\dot{B}_{r,\infty}^{-\varepsilon+1}}\|\theta\|_{B_{\infty,1}^{\frac{2}{r}+2\varepsilon-1,g}}+\|u\|_{L^2}\|\theta\|_{L^2}\Big)\mathrm{d}\tau \\
   &\leq C\int_{0}^{t}\Big\{(\|u\|_{L^2}+\|\omega\|_{B_{r,\infty}^{-\varepsilon}})\|\theta\|_{L^\infty}+\|u\|_{L^2}\|\theta\|_{L^2}\Big\}\mathrm{d}\tau \\
   &\leq C(T,u_{0},\theta_{0}).
 \end{split}
 \end{equation}
 We next estimate the last term in the right hand side of ({\ref{3.22}}). Invoking the Bony decomposition, the term $\|[\Delta_{j},u \cdot \nabla]G\|_{L^r}$ can be written as
      \begin{equation}\label{3.25}
      \begin{split}
 \|[\Delta_{j},u \cdot \nabla]G\|_{L^r} &= \sum_{|k-j| \leq 4} [\Delta_{j}, S_{k-1}u \cdot \nabla]\Delta_{k}G+\sum_{|k-j| \leq 4} [\Delta_{j}, \Delta_{k}u \cdot \nabla]S_{k-1}G \\
 &+ \sum_{k \geq j-4}[\Delta_{j},\Delta_{k}u \cdot \nabla]\tilde{\Delta}_{k} G \\
 &\triangleq I_{1}+I_{2}+I_{3}.
  \end{split}
  \end{equation}
  By Lemma {\ref{lem 2.5}}, one has
        \begin{equation}\label{3.26}
        \begin{split}
        \|I_{1}\|_{L^r} &\leq C\||x|^{1-\varepsilon} h_{j}\|_{L^1}\|S_{j-1}u\|_{\dot{B}_{r,\infty}^{1-\varepsilon}}\|\nabla \Delta_{j}G\|_{L^\infty}  \\
         &\leq C2^{j(\varepsilon+\frac{2}{r})}\|\omega\|_{\dot{B}_{r,\infty}^{-\varepsilon}}\|\Delta_{j}G\|_{L^r},
  \end{split}
  \end{equation}
        \begin{equation}\label{3.27}
        \begin{split}
        \|I_{2}\|_{L^r}  &\leq C2^{- j}\|\nabla S_{j-1}G\|_{L^\infty} \|\Delta_{j}(\nabla u)\|_{L^r} \\
        &\leq C2^{-j}2^{\varepsilon j}\|\omega\|_{\dot{B}_{r,\infty}^{-\varepsilon}}\sum_{m \leq j-2}2^{m}2^{\frac{2}{r}m}\|\Delta_{m} G\|_{L^r}  \\
        &\leq C2^{j(\varepsilon+\frac{2}{r})}\|\omega\|_{\dot{B}_{r,\infty}^{-\varepsilon}}\sum_{m\leq j-2} 2^{(m-j)(\frac{2}{r}+1)}\|\Delta_{m}G\|_{L^r}.
  \end{split}
  \end{equation}
  We further split the term $I_{3}$ into two parts
    \begin{align}
 I_{3} \triangleq I_{31}+I_{32},\nn
  \end{align}
  where
    \begin{align}
 I_{31}= \sum_{k \geq j-4,k \geq 0}[\Delta_{j},\Delta_{k}u \cdot \nabla]\tilde{\Delta}_{k} G,\quad I_{32}= [\Delta_{j},\Delta_{-1}u \cdot \nabla] \tilde{\Delta}_{-1} G.\nn
  \end{align}
The aim of splitting of $I_3$ is such that $I_{31}$ removes the low-frequency term $\Delta_{-1}$ and $I_{32}$ includes the low-frequency term $\Delta_{-1}$. By Bernstein's inequality,
 \begin{equation}\label{3.28}
 \begin{split}
 \|I_{31}\|_{L^r} &= \Big\|\sum_{k-j \geq -4,k \geq 0} \Big\{\Delta_{j} \nabla \cdot (\Delta_{k}u \tilde{\Delta}_{k} G)-\Delta_{k}u \Delta_{j} \nabla \cdot \tilde{\Delta}_{k} G \Big\}\Big\|_{L^r}  \\
  &\leq C\sum_{k-j \geq -4,k \geq 0}\Big(2^{j}\|\Delta_{k}u\tilde{\Delta}_{k}G\|_{L^r}+\|\Delta_{k}u\|_{L^r}\|\Delta_{j}\nabla \tilde{\Delta}_{k}G\|_{L^\infty}\Big)\\
   &\leq C\sum_{k-j \geq -4,k \geq 0}2^{j}\|\Delta_{k}u\|_{L^r}\|\Delta_{k}G\|_{L^\infty} \\
   &\leq C\sum_{k-j \geq -4 ,k \geq 0}2^{j}2^{(\varepsilon-1)k}2^{\frac{2}{r}k}\|\Delta_{k}\Lambda^{1-\varepsilon}u\|_{L^r}\|\Delta_{k}G\|_{L^r}\\
   &\leq C2^{j(\varepsilon+\frac{2}{r})}\|\omega\|_{\dot{B}_{r,\infty}^{-\varepsilon}}\sum_{k-j \geq 4, k \geq 0}2^{(j-k)(1-\frac{2}{r}-\varepsilon)}\|\Delta_{k}G\|_{L^r}.
   \end{split}
    \end{equation}
 For $I_{32}$, by ({\ref{2.3}}), we derive
 \begin{equation}\label{3.29}
 \begin{split}
 \|I_{32}\|_{L^r} &\leq C\||x|h_j\|_{L^1}\|\Delta_{-1}\nabla u\|_{L^r}\|\tilde{\Delta}_{-1} \nabla G\|_{L^\infty} \\
    &\leq C2^{-j}\|u\|_{L^2}\|G\|_{L^2}.
    \end{split}
    \end{equation}
 Combining ({\ref{3.26}}), ({\ref{3.27}}), ({\ref{3.28}}) and ({\ref{3.29}}) into ({\ref{3.25}}), we have
     \begin{equation}\label{3.30}
     \begin{split}
    C\sum_{j \geq -1}2^{j(\frac{2}{r}+\varepsilon-1)}\|[\Delta_{j},u \cdot \nabla]G\|_{L^r}
    &\leq C\sum_{j \geq -1}2^{j(\frac{2}{r}+2\varepsilon-1)}\|\omega\|_{\dot{B}_{r,\infty}^{-\varepsilon}}2^{j\frac{2}{r}}\Big\{\|\Delta_{j}G\|_{L^r}\\&+\sum_{m\leq j-2} 2^{(m-j)(\frac{2}{r}+1)}\|\Delta_{m}G\|_{L^r} + \sum_{k \geq j -4}2^{(j-k)(1-\frac{2}{r}-\varepsilon)}\|\Delta_{k}G\|_{L^r}\Big\}\\&+C\|u\|_{L^2}\|G\|_{L^2},
    \end{split}
  \end{equation}
where we choose $\varepsilon$  small enough such that
\begin{equation}\label{negative}
\frac{2}{r}+2\varepsilon-1 <0.
\end{equation}
Inserting ({\ref{3.23}}), ({\ref{3.24}}) and ({\ref{3.30}}) into ({\ref{3.22}}), we have
\begin{equation}\label{3.31}
\begin{split}
\|G\|_{L_{t}^{1}B_{r,1}^{\frac{2}{r}}} &\leq C(T) + C\sum_{j \geq -1}2^{j(\frac{2}{r}+2\varepsilon-1)}\|\omega\|_{L_{t}^{1}\dot{B}_{r,\infty}^{-\varepsilon}}2^{j\frac{2}{r}}\Big\{\|\Delta_{j}G\|_{L_{t}^{1}L^r}\\&+\sum_{m\leq j-2} 2^{(m-j)(\frac{2}{r}+1)}\|\Delta_{m}G\|_{L_{t}^{1}L^r} + \sum_{k \geq j -4}2^{(j-k)(1-\frac{2}{r}-\varepsilon)}\|\Delta_{k}G\|_{L_{t}^{1}L^r}\Big\}.
\end{split}
\end{equation}
We decompose  the right hand side of ({\ref{3.31}}) into two parts K and L as follows
\begin{align}
&K= C(T) + C\sum_{j = -1}^{N}2^{j(\frac{2}{r}+2\varepsilon-1)}\|\omega\|_{L_{t}^{1}\dot{B}_{r,\infty}^{-\varepsilon}}2^{j\frac{2}{r}}\Big\{\|\Delta_{j}G\|_{L_{t}^{1}L^r} \nn \\ &+\sum_{m\leq j-2} 2^{(m-j)(\frac{2}{r}+1)}\|\Delta_{m}G\|_{L_{t}^{1}L^r} + \sum_{k \geq j -4}2^{(j-k)(1-\frac{2}{r}-\varepsilon)}\|\Delta_{k}G\|_{L_{t}^{1}L^r}\Big\}, \nn\\
&L= C\sum_{j > N}2^{j(\frac{2}{r}+2\varepsilon-1)}\|\omega\|_{L_{t}^{1}\dot{B}_{r,\infty}^{-\varepsilon}}\Big\{2^{j\frac{2}{r}}\|\Delta_{j}G\|_{L_{t}^{1}L^r} \nn \\ &+\sum_{m\leq j-2} 2^{(m-j)}2^{m\frac{2}{r}}\|\Delta_{m}G\|_{L_{t}^{1}L^r} + \sum_{k \geq j -4}2^{(j-k)(1-\varepsilon)}2^{k\frac{2}{r}}\|\Delta_{k}G\|_{L_{t}^{1}L^r}\Big\}, \nn
\end{align}
where the integer $N$ will be chosen later.
Noting the formula (\ref{negative}), by the discrete Young's inequality, the term $K$ can be estimated as follows
\begin{equation}\label{3.32}
\begin{split}
K &\leq C(T)+C2^{\frac{2}{r}N}\|\omega\|_{L_{t}^{1}\dot{B}_{r,\infty}^{-\varepsilon}}\|G\|_{L_{t}^{1}L^r}\\
&\leq C(T)+C2^{\frac{2}{r}N}\|\omega\|_{L_{t}^{1}\dot{B}_{r,\infty}^{-\varepsilon}}.
\end{split}
\end{equation}
Similarly the term $L$ can be estimated as follows
\begin{equation}\label{3.33}
L \leq  C2^{N(\frac{2}{r} + 2\varepsilon -1)}\|\omega\|_{L_{t}^{1}\dot{B}_{r,\infty}^{-\varepsilon}}\|G\|_{L_{t}^{1}B_{r,1}^{\frac{2}{r}}}.
\end{equation}
Now take $N$ large enough such that
\begin{align}
C2^{N(\frac{2}{r} + 2\varepsilon -1)}\|\omega\|_{L_{t}^{1}\dot{B}_{r,\infty}^{-\varepsilon}} \leq \frac{1}{2}.\nn
\end{align}
Inserting (\ref{3.32}) and (\ref{3.33}) into (\ref{3.31}), we get
\begin{align}
\|G\|_{L_{t}^{1}B_{r,1}^{\frac{2}{r}}} &\leq C(T) + C2^{\frac{2}{r}N}\|\omega\|_{L_{t}^{1}\dot{B}_{r,\infty}^{-\varepsilon}} \nn\\
&\leq C(T,u_{0},\theta_{0}),\nn
\end{align}
which completes the proof of Proposition {\ref{pro 3.4}}.
\end{proof}

\subsection{The estimates of $\|\omega\|_{L^{1}_{t}{B_{\infty,1}^{0}}}$, $\|\nabla u\|_{L_{t}^{1}L^\infty}$ and $\|\nabla \theta\|_{L^\infty}$}
      In this subsection first shows that $\|\omega\|_{L^{1}_{t}{B_{\infty,1}^{0}}}$, and further obtains the estimate of $\|\nabla u\|_{L_{t}^{1}L^\infty}$, finally the estimate of $\|\nabla \theta\|_{L^\infty}$ is deduced.
\begin{proposition}\label{pro 3.5}
Assume $(u_{0},\theta_{0})$ satisfies the assumptions in Theorem {\ref{thm 1.1}}.  Then we have the following global $\textit{a priori}$ bounds,
 \begin{align}
    \|\omega\|_{L^{1}_{t}{B_{\infty,1}^{0}}} \leq C(T,u_{0},\theta_{0}), \quad \|\nabla u\|_{L_{t}^{1}{L^\infty}}\leq C(T,u_{0},\theta_{0}), \quad \|\nabla \theta\|_{L^\infty} \leq C(T,u_{0},\theta_{0}). \nn
     \end{align}
\end{proposition}
\begin{proof}
First, we prove the bound of $ \|\omega\|_{L^{1}_{t}{B_{\infty,1}^{0}}}$. Taking $r \in (2,4)$, we have
\begin{equation}\label{omega}
 \begin{split}
\|\omega\|_{L^{1}_{t}{B_{\infty,1}^{0}}} &\leq \|G\|_{L^{1}_{t}{B_{\infty,1}^{0}}} + \|\mathcal{R}_{g}\theta\|_{L^{1}_{t}{B_{\infty,1}^{0}}}\\
&\leq \sum_{j \geq -1}2^{\frac{2}{r}j}\|\Delta_{j}G\|_{L^{1}_{t}{L^r}}+C\|\theta\|_{L^{1}_{t}{B_{\infty,1}^{0,g}}}\\
&\leq \|G\|_{L^{1}_{t}{B_{r,1}^{\frac{2}{r}}}}+C\|\theta\|_{L^{1}_{t}{B_{\infty,1}^{0,g}}}.
 \end{split}
\end{equation}
By Lemma {\ref{lem 2.9}} and Proposition \ref{3.3}, we deduce that
 \begin{equation} \label{3.34}
 \begin{split}
\|\theta\|_{B_{\infty,1}^{0,g}} &\leq C(\|\theta_{0}\|_{B_{\infty,1}^{0,g}}+\|u\|_{L^{1}_{t}{B_{\infty,1}^{0,g}}})\Big(1+\int_{0}^{t}\|\nabla u\|_{L^\infty}\mathrm{d}\tau\Big)\\
&\leq C(\|\theta_{0}\|_{B_{\infty,1}^{0,g}}+\|u\|_{L^{1}_{t}{L^2}}+\|\omega\|_{L^{1}_{t}{B_{\infty,1}^{\varepsilon-1}}})(1+\|u\|_{L^{1}_{t}{L^{2}}}+\|\omega\|_{L^{1}_{t}{B_{\infty,1}^{0}}}) \\
&\leq C(\|\theta_{0}\|_{B_{\infty,1}^{0,g}}+\|u\|_{L^{1}_{t}{L^2}}+\|G\|_{L^{1}_{t}{L^r}}+\|\theta\|_{L^{1}_{t}{L^r}})(1+\|u\|_{L^{1}_{t}{L^{2}}}+\|\omega\|_{L^{1}_{t}{B_{\infty,1}^{0}}}) \\
&\leq C(\|\theta_{0}\|_{B_{\infty,1}^{0,g}}+1)(1+\|\omega\|_{L^{1}_{t}{B_{\infty,1}^{0}}}).
     \end{split}
     \end{equation}
Inserting the estimate (\ref{3.34}) into (\ref{omega}), by Proposition \ref{pro 3.4}, we have
 \begin{align}
\|\omega\|_{L^{1}_{t}{B_{\infty,1}^{0}}} &\leq C+C(\|\theta_{0}\|_{B_{\infty,1}^{0,g}}+1)T  + C(\|\theta_{0}\|_{B_{\infty,1}^{0,g}}+1)\int_{0}^{t}\|\omega\|_{L^{1}_{\tau}{B_{\infty,1}^{0}}}\mathrm{d}\tau.\nn
\end{align}
By Gr$\rm\ddot{o}$nwall's inequality, we have
\begin{align}
\|\omega\|_{L^{1}_{t}{B_{\infty,1}^{0}}} \leq C(T, u_{0}, \theta_{0}).\nn
\end{align}
 Clearly one has
$$
\|\nabla u\|_{L_{t}^{1}{L^\infty}} \le \|u\|_{L^{1}_{t}{L^{2}}}+\|\omega\|_{L^{1}_{t}{B_{\infty,1}^{0}}} \leq  C(T, u_{0}, \theta_{0}).
$$
Next, we estimate $\|\nabla\theta\|_{L^\infty}$. Applying the operator $\nabla$ on the both sides of the equation $({\ref{1.1}})_{2}$ and dotting it by $|\nabla \theta|^{r-2}\nabla \theta$ and integrating on $\mathbb{R}^{2}$ for any $2<r<\infty$, we have
\begin{align}
\frac{1}{r}\frac{\mathrm{d}}{\mathrm{d}t}\|\nabla \theta\|_{L^r}^{r} &\leq \int_{\mathbb{R}^{2}} \nabla u_{2}
\cdot |\nabla \theta|^{r-2} \nabla \theta \mathrm{d}x -\int_{\mathbb{R}^{2}}\nabla u \cdot \nabla \theta \cdot |\nabla \theta|^{r-2} \nabla \theta \mathrm{d}x \nn\\
&\leq \|\nabla u\|_{L^r}\|\nabla \theta\|_{L^r}^{r-1} + \|\nabla u\|_{L^\infty}\|\nabla \theta\|_{L^r}^{r}.\nn
\end{align}
Thus, we have
\begin{align}
\frac{\mathrm{d}}{\mathrm{d}t}\|\nabla \theta\|_{L^r} \leq \|\nabla u\|_{L^r} + \|\nabla u\|_{L^\infty}\|\nabla \theta\|_{L^r}.\nn
\end{align}
Applying Gr$\rm\ddot{o}$nwall's inequality it yields
\begin{align}
\|\nabla \theta\|_{L^r} &\leq (\|\nabla \theta_{0}\|_{L^r} + \int_{0}^{t} \|\nabla u\|_{L^r}\mathrm{d}\tau)\exp\Big\{\int_{0}^{t}\|\nabla u\|_{L^\infty}\mathrm{d}\tau\Big\} \leq C(T,u_{0}, \theta_{0}).\nn
\end{align}
Let $r \rightarrow \infty$, we have
\begin{align}
\|\nabla \theta\|_{L^\infty}  \leq C(T,u_{0}, \theta_{0}),\nn
\end{align}
which completes the proof of Proposition {\ref{pro 3.5}}.
\end{proof}

\subsection{The global $H^{s}$ estimates of $u$ and $\theta$}
     In this subsection, we establish the global $H^{s}$ estimates of $u$ and $\theta$ for $s\ge 2$.
\begin{proposition}\label{pro 3.6}
Assume $(u_{0},\theta_{0})$ satisfies the assumptions in Theorem {\ref{thm 1.1}}. Then the corresponding solution $(u,\theta)$ of ({\ref{1.1}}) is bounded in $H^{s}(\mathbb{R}^{2})$.
\end{proposition}
\begin{proof}
We apply $\Lambda^{s}$ to the equations $u$ and $\theta$, and take the $L^2$ inner product of the resulting equations with $(\Lambda^{s}u,\Lambda^{s}\theta)$ to obtain
\begin{equation}\label{3.37}
\begin{split}
&\frac{1}{2}\frac{\mathrm{d}}{\mathrm{d}t}(\|\Lambda^{s}u\|_{L^2}^{2}+\|\Lambda^{s}\theta\|_{L^2}^{2})+\|\mathcal{L}^{\frac{1}{2}}\Lambda^{s}u\|_{L^2}^{2} \\
 &\leq -\int_{\mathbb{R}^{2}}[\Lambda^{s},u \cdot \nabla]u \cdot \Lambda^{s}u \mathrm{d}x+\int_{\mathbb{R}^{2}}\Lambda^{s}(\theta \cdot e_{2}) \cdot \Lambda^{s}u \mathrm{d}x \\
 &\ -\int_{\mathbb{R}^{2}}[\Lambda^{s},u \cdot \nabla]\theta \cdot \Lambda^{s} \theta \mathrm{d}x+\int_{\mathbb{R}^{2}}\Lambda^{s}u_{2} \cdot \Lambda^{s}\theta \mathrm{d}x \\
 &:=A_{1}+A_{2}+A_{3}+A_{4}.
 \end{split}
\end{equation}
By Lemma {\ref{lem 2.6}}, we can deduce that
\begin{equation*}
\begin{split}
|A_{1}| &\leq \|[\Lambda^{s},u \cdot \nabla]u\|_{L^2}\|\Lambda^{s}u\|_{L^2} \\
&\leq C\|\nabla u\|_{L^\infty}\|\Lambda^{s}u\|_{L^2}^{2},\\
|A_2+&A_4| \leq (\|\Lambda^{s}\theta\|_{L^2}^{2} + \|\Lambda^{s}u\|_{L^2}^{2}),
\end{split}
\end{equation*}
\begin{equation*}
\begin{split}
|A_{3}| &\leq \|[\Lambda^{s},u \cdot \nabla]\theta\|_{L^2}\|\Lambda^{s}\theta\|_{L^2} \\
&\leq C(\|\nabla u\|_{L^\infty}\|\Lambda^{s}\theta\|_{L^2}+\|\Lambda^{s}u\|_{L^2}\|\nabla \theta\|_{L^\infty})\|\Lambda^{s}\theta\|_{L^2}\\
&\leq C(\|\nabla u\|_{L^\infty}+1)(\|\Lambda^{s}u\|_{L^2}^{2}+\|\Lambda^{s}\theta\|_{L^2}^{2}),
\end{split}
\end{equation*}
Inserting the above estimates into ({\ref{3.37}}), we have
\begin{equation*}
\begin{split}
&\frac{\mathrm{d}}{\mathrm{d}t}(\|\Lambda^{s}u\|_{L^2}^{2}+\|\Lambda^{s}\theta\|_{L^2}^{2})+2\|\mathcal{L}^{\frac{1}{2}}\Lambda^{s}u\|_{L^2}^{2} \\
 &\leq C(\|\nabla u\|_{L^\infty}+1)(\|\Lambda^{s}u\|_{L^2}^{2}+\|\Lambda^{s}\theta\|_{L^2}^{2}).
 \end{split}
 \end{equation*}
 Applying Gr$\rm\ddot{o}$nwall's inequality it yields
 \begin{equation*}
\begin{split}
&(\|u\|_{H^s}^{2}+\|\theta\|_{H^s}^{2})+2\int_{0}^{t}\|\mathcal{L}^{\frac{1}{2}}u\|_{H^s}^{2}\mathrm{d}\tau \\
 &\leq (\|u_{0}\|_{H^s}^{2}+\|\theta_{0}\|_{H^s}^{2}) \exp\{\int_{0}^{t}(\|\nabla u\|_{L^\infty}+1)\mathrm{d}\tau \}\\
 &\leq C(T,u_{0},\theta_{0}).
 \end{split}
 \end{equation*}
The proof of Proposition {\ref{pro 3.6}} is thus completed.
\end{proof}

\section{Proof of Theorem 1.1}
\setcounter{equation}{0}

In this section we only sketch the existence proof of Theorem {\ref{thm 1.1}} because it is a standard process, and then we prove the uniqueness of solutions concisely.

\begin{proof}
First, we prove the existence of a solution. Let $N > 0$ be an integer, we define $J_N$ as
\begin{align}
\mathcal{F}({J_{N}f}) = \chi_{B(0,N)} (\xi)\mathcal{F}(f)(\xi),\nn
\end{align}
where
$\chi_{B(0,N)}(\xi) =
\begin{cases}
1, \mbox {if} \  \xi \in B(0,N),\\
0, \mbox{if}\ \xi  \notin B(0,N).
\end{cases}$\par
We consider an approximate equation to the equation ({\ref{1.1}})
\begin{eqnarray}\label{4.1}
	\left\{\begin {array}{lll}
	\partial_{t}u^{N}+\mathbb{P}J_{N}(J_{N}u^{N} \cdot \nabla (J_{N}u^{N}))+J_{N}\mathcal{L} u^{N}=J_{N}\theta^{N}e_{2},\\
	\partial_{t}\theta^{N}+J_{N}(J_{N}u^{N} \cdot \nabla J_{N}\theta^{N})=J_{N}u^{N}_{2},\\
	\nabla \cdot u^{N}=0,\\
	(u^{N},\theta^{N})(x,t)_{t=0} = (J_{N}u_{0},J_{N}\theta_{0}),
	\end{array}\right.
	\end{eqnarray}
where $\mathbb{P}$ denotes the Leray projection operator (onto divergence-free vector fields).  By the compactness argument, this system of equations admits a global solution $(u_{N}, \theta_{N})$ that satisfies the global \textit{a priori} estimates stated in Proposition {\ref{pro 3.6}} uniformly in terms of $N$. By making $N$ approach $\infty$ and going through a limit process, we conclude that $(u_{N}, \theta_{N})$ converges to a solution of \eqref{1.1}. Since these processes are standard, we omit further details, the specific processes can be found in Majda and Bertozzi \cite{MB2001}. The global existence of a solution is thus obtained.\par
 Next, we show the uniqueness of the solution. Assume that $(u^{(1)},\theta^{(1)})$ and
$(u^{(2)},\theta^{(2)})$ are two solutions of the equations \eqref{1.1} emanating from the same initial data $(u_0,\theta_0)$. Their difference
$(\bar{u},\bar{\theta})$ denoted by
$$\bar{u}=u^{(2)}-u^{(1)},\ \bar{p}=p^{(2)}-p^{(1)},\ \bar{\theta}=\theta^{(2)}-\theta^{(1)}$$
satisfy the following equations
   \begin{equation}\label{4.2}
\left\{\aligned
&\partial_{t}\bar{u}+u^{(1)} \cdot \nabla \bar{u}+\bar{u} \cdot \nabla u^{(2)}+\mathcal{L}\bar{u}=-\nabla \bar{p}+\bar{\theta} e_{2},\\
&\partial_{t} \bar{\theta}+u^{(1)} \cdot \nabla \bar{\theta}+\bar{u} \cdot \nabla \theta^{(2)}=\bar{u} \cdot e_{2}.
\endaligned\right.
\end{equation}
Dotting ({\ref{4.2}}) by $(u,\theta)$, we deduce
 \begin{equation}\label{4.3}
\frac{1}{2}\frac{\mathrm{d}}{\mathrm{d}t}(\|\bar{u}\|_{L^2}^{2}+\|\bar{\theta}\|_{L^2}^{2})+\|\mathcal{L}^{\frac{1}{2}}\bar{u}\|_{L^2}^{2} = M_{1} + M_{2} + M_{3} + M_{4},
\end{equation}
where
\begin{align*}
&M_{1}=-\int_{\mathbb{R}^{2}}\bar{u} \cdot\nabla u^{(2)}\cdot \bar{u} \mathrm{d}x,\qquad
M_{2}=\int_{\mathbb{R}^{2}}\bar{\theta} e_{2} \cdot \bar{u} \mathrm{d}x,\\
&M_{3}=-\int_{\mathbb{R}^{2}}\bar{u} \cdot \nabla \theta^{(2)} \cdot \bar{\theta} \mathrm{d}x,\qquad
M_{4}=\int_{\mathbb{R}^{2}}\bar{u} \cdot \bar{\theta} \mathrm{d}x.
\end{align*}
By the H\"{o}lder and Young's inequalities, one has
\begin{align*}
|M_{1}| &\leq \|\nabla u^{(2)}\|_{L^{\infty}}\|\bar{u}\|_{L^{2}}^{2},
\end{align*}
\begin{align*}
|M_{2}|+|M_{4}| &\leq C(\|\bar{u}\|_{L^{2}}^{2}+\|\bar{\theta}\|_{L^{2}}^{2}),
\end{align*}
\begin{align*}
|M_{3}| \leq C(\|\bar{u}\|_{L^2}^{2}+\|\bar{\theta}\|_{L^2}^{2})\|\nabla \theta^{(2)}\|_{L^\infty}.
\end{align*}
Inserting above estimates into (\ref{4.3}), we have
\begin{align*}
\frac{\mathrm{d}}{\mathrm{d}t}(\|\bar{u}\|_{L^2}^{2}+\|\bar{\theta}\|_{L^2}^{2})+2\|\mathcal{L}^{\frac{1}{2}}\bar{u}\|_{L^2}^{2} \leq C(1+\|\nabla u^{(2)}\|_{L^\infty}+\|\nabla \theta^{(2)}\|_{L^\infty})(\|\bar{u}\|_{L^2}^{2}+\|\bar{\theta}\|_{L^2}^{2}).
\end{align*}
By the the boundedness  of $\int^T_0\Big(1+\|\nabla u^{(2)}\|_{L^\infty}+\|\nabla \theta^{(2)}\|_{L^\infty}\Big)\mbox{d}t$ and $u_0=b_0\equiv 0$, applying Gr$\rm\ddot{o}$nwall's inequality it yields
\begin{align*}
\|\bar{u}\|_{L^2}^{2}+\|\bar{\theta}\|_{L^2}^{2}  \leq C(\|\bar{u}_{0}\|_{L^2}^{2}+\|\bar{\theta}_{0}\|_{L^2}^{2}) \exp \int^T_0\Big(1+\|\nabla u^{(2)}\|_{L^\infty}+\|\nabla \theta^{(2)}\|_{L^\infty}\Big)\mbox{d}t=0,
\end{align*}
and we thus complete the proof of Theorem {\ref{thm 1.1}}.
\end{proof}
\vspace{.1in}

\textbf{Declarations}
\vspace{.1in}

\textbf{Data Availability Statement} Data sharing not applicable to this article as no datasets were generated or analysed during the current study.

\textbf{Conflict of interest} On behalf of all authors, the corresponding author states that there
is no conflict of interest.


\begin{thebibliography}{9}
         \bibitem{DCJ2011} D. Adhikari, C. Cao and J. Wu. Global regularity results for the 2D Boussinesq equations with vertical dissipation.
     \textit{J. Differential Equations,} \textbf{251} (2011), 1637-1655.

             \bibitem{HJR2011}H. Bahouri, Y. Chemin and R. Danchin. Fourier Analysis and Nonlinear Partial Differential Equations. In: Grundlehren der Mathematischen Wissenschaften, Vol. 343, \textit{Springer,} 2011.

     \bibitem{J1998}Y. Chemin. Perfect Incompressible Fluids.
 \textit{Oxford University Press,} 1998.

         \bibitem{QCZ2007}L. Chen, C. Miao and Z. Zhang. A new Bernstein inequality and the 2D dissipative quasigeostrophic equation.
      \textit{Comm. Math. Phys.} \textbf{ 271} (2007), 821-838.

           \bibitem{DJ2012}D. Chae, J. Wu. The 2D Boussinesq equations with logarithmically supercritical velocities.
\textit{Advances in Math.} \textbf{230} (2012), 1618-1645.

  \bibitem{RM2011}R. Danchin, M. Paicu. Global existence results for the anisotropic Boussinesq system in dimension two.
     \textit{Math. Models Methods Appl. Sci.} \textbf{21} (2011), 421-457.

      \bibitem{MALV2014}M. Dabkowski, A. Kiselev, L. Silvestre and V. Vicol. Global well-posedness of slightly supercritical active scalar equations.
      \textit{ Analysis and PDE} \textbf{7} (2014),43-72


           \bibitem{COR87}C. Foias, O. Manley and R. Temam. Attractors for the B\'{e}nard problem: existence and physical bounds on their fractal dimension.
      \textit{Nonlinear Anal.} \textbf{11(8)} (1987), 939–967.

              \bibitem{AME2015}  A. Farhat, M. Jolly and E. Titi. Continuous data assimilation for the 2D B\'{e}nard convection through velocity measurements alone.
     \textit{Physica D.} \textbf{303} (2015), 59–66.


     \bibitem{TSF2010}T. Hmidi, S. Keraani and F. Rousset. Global well-posedness for a Boussinesq-Navier-Stokes system with critical dissipation.
      \textit{J. Differential Equations,} \textbf{249} (2010),  2147-2174.

       \bibitem{TSF2011}T. Hmidi, S. Keraani and F. Rousset. Global well-posedness for Euler-Boussinesq system with critical dissipation.
      \textit{ Comm. Partial Differential Equations,} \textbf{36} (2011), 420-445.


     \bibitem{AH2020}A. Hanachi, H. Houamed and M. Zerguine. On the global well-posedness of the axisymmetric viscous Boussinesq system in critical Lebesgue spaces.
     \textit{Discrete Contin. Dyn. Syst.} \textbf{40(11)} (2020), 6473-6506.

           \bibitem{QCJZ2014}Q. Jiu, C. Miao, J. Wu and Z. Zhang. The 2D incompressible Boussinesq equations with general critical dissipation.
\textit{ SIAM J. Math. Anal.} \textbf{46} (2014), 3426-3454.

    \bibitem{KP1988} T. Kato, G. Ponce. Commutator estimates and the Euler and Navier-Stokes equations.
\textit{Comm. Pure Appl. Math.} \textbf{41(7)} (1988), 891-907.

          \bibitem{KPV1991} C. Kenig, G. Ponce and L. Vega. Well-posedness of the initial value problem for the Korteweg-de Vries equations.
    \textit{J. Amer. Math. Soc.} \textbf{4(2)} (1991), 323-347.

     \bibitem{MB2001} A. Majda, A. Bertozzi. Vorticity and Incompressible Flow.
    \textit{Cambridge University Press,} Cambridge, 2001.

         \bibitem{TS07}T. Ma, S. Wang. Rayleigh B\'{e}nard convection: dynamics and structure in the physical space.
     \textit{ Commun. Math. Sci.} \textbf{5 (3)} (2007), 553-574.

       \bibitem{CL2011} C. Miao, L. Xue. On the global well-posedness of a class of Boussinesq- Navier-Stokes systems.
            \textit{NoDEA  Nonlinear Differential Equations Appl.} \textbf{18} (2011), 707-735.


            \bibitem{DLJD}D. KC, D. Regmi, L. Tao, and J. Wu. The 2D Boussinesq-Navier-Stokes equations with logarithmically supercritical dissipation.  \textit{J. Math. Study.} \textbf{57} (2024), 101-132.

             \bibitem{P1968}H. Rabinowitz. Existence and nonuniqueness of rectangular solutions of the B\'{e}nard problem.
      \textit{Arch. Ration. Mech. Anal.} \textbf{29} (1968), 32-57.

     \bibitem{Z2014}Z. Ye. Blow-up criterion of smooth solutions for the Boussinesq equations.
     \textit{Nonlinear Anal.} \textbf{110} (2014), 97-103.

               \bibitem{Z2017}Z. Ye. Regularity criterion of the 2D B\'{e}nard equations with critical and supercritical disspipation.
     \textit{Nonlinear Anal.} \textbf{156} (2017), 111-143.
\end{thebibliography}
\end{document}